\documentclass[12pt, reqno]{amsart}
\usepackage{amsmath, amsthm, amscd, amsfonts, amssymb, graphicx, color}
\usepackage[bookmarksnumbered, colorlinks, plainpages]{hyperref}

\usepackage[all]{xy}
\newtheorem{theorem}{Theorem}[section]

\newtheorem{proposition}[theorem]{Proposition}
\newtheorem{corollary}[theorem]{Corollary}
\theoremstyle{definition}
\newtheorem{definition}[theorem]{Definition}
\newtheorem{example}[theorem]{Example}

\theoremstyle{remark}
\newtheorem{remark}[theorem]{Remark}
\numberwithin{equation}{section}

\begin{document}
	\setcounter{page}{1}
	\title[\textit{DW}-DP operators and \textit{DW}-limited operators]{\textit{DW}-DP operators and \textit{DW}-limited operators on Banach lattices}
	
    \author[J.X. Chen]
       {Jin Xi Chen}
       \address{School of Mathematics, Southwest Jiaotong University, Chengdu 610031, China}
       \email{jinxichen@swjtu.edu.cn}

	 \author[J. Feng]
	  {Jingge Feng}
	 \address { School of Mathematics, Southwest Jiaotong University, Chengdu 610031, China}
	
\email{fjg0825@my.swjtu.edu.cn}

\subjclass[2010] {Primary 46B42; Secondary 46B50, 47B65}

	\keywords{disjointly weakly compact set, $DW$-DP operator, \textit{DW}-limited operator, weak Dunford-Pettis operator, weak$^*$ Dunford-Pettis operator, Banach lattice.}
	
\begin{abstract}
This paper is devoted to the study of  two classes of operators related to disjointly weakly compact sets,  which we call $DW$-DP operators and $DW$-limited operators, respectively. They carry disjointly weakly compact subsets of a Banach lattice onto Dunford-Pettis sets and limited sets, respectively. We show that $DW$-DP (resp. $DW$-limited) operators are precisely the operators which are both weak Dunford-Pettis  and order Dunford-Pettis (resp.  weak$^*$ Dunford-Pettis and  order limited) operators. Furthermore, the approximation properties of positive $DW$-DP and positive $DW$-limited operators are given.
\end{abstract}
\maketitle\baselineskip 4.80mm
	
\section{Introduction and preliminaries}
Throughout this paper, $X, Y$ will denote real Banach spaces, and
	$E, F$ will denote real Banach lattices. $B_X$ is the closed unit ball of $X$, and $E^+$ is the positive cone of the Banach lattice $E$.  The solid hull of a subset $A$ of $E$ is denoted by $Sol(A):=\lbrace y\in E:|y| \leq |x|$ \ for \ some\ $x\in A\rbrace$.

Recall that a bounded subset $A$ of  $X$ is called a \textit{Dunford-Pettis} (resp. \textit{limited}) \textit{set} if every weakly null (resp. weak$^*$-null) sequence $(x_n')$ of $X^\prime$ converges uniformly to zero on $A$, i.e., $\sup_{x\in{A}}{\vert{x_n'(x)}\vert}\to{0}\,\, (n\to\infty)$, or equivalently, if  every weakly compact (resp. continuous) operator $T:X\rightarrow c_0$ carries $A$ to a relatively compact set (see \cite{positiveoperators,DPset,Limitedsets}). $X$ is said to have the \textit{Dunford-Pettis property} (resp. \textit{DP$^{\,*}$ property}) if every relatively weakly compact subset of $X$ is a Dunford-Pettis  (resp. limited) set, or equivalently, if each weakly compact  (resp. continuous) linear operator from $X$ to $c_0$ is a Dunford-Pettis operator \cite{DPandDPstarproperty, ocDPstar}. Note that a Dunford-Pettis operator  maps relatively weakly compact sets onto relatively compact ones. Clearly, a relatively compact subset of a Banach space is a Dunford-Pettis and limited set. Therefore, Aliprantis and Burkinshaw \cite{wDPoperator}, and  El Kaddouri et al. \cite{KHBM, w*DPoperator} introduced the class of operators which map relatively weakly compact sets  onto  Dunford-Pettis (resp. limited) sets . A bounded linear operator $T:X\to{Y}$ between Banach spaces is called
 \begin{itemize}
   \item[$\raisebox{0mm}{------}$]  a \textit{weak Dunford-Pettis} (abbr. wDP) \textit{operator}  if $x_n\xrightarrow{w}0$ in $X$ and $y_n'\xrightarrow{w}0$ in $Y'$ imply $y_n'(Tx_n)\rightarrow0$, or equivalently, if $T$ carries relatively weakly compact subsets of $X$ onto  Dunford-Pettis subsets of $Y$ (see \cite{wDPoperator, positiveoperators}).
   \item[$\raisebox{0mm}{------}$]  a \textit{weak$^{\,*}$ Dunford-Pettis} (abbr. w$^*$DP) \textit{operator} whenever $x_n\xrightarrow{w}0$ in $X$ and $y_n'\xrightarrow{w^*}0$ in $Y'$ imply $y_n'(Tx_n)\rightarrow0$, or equivalently, whenever $T$ carries  relatively weakly compact subsets of $X$ onto  limited subsets of $Y$ \cite{KHBM, w*DPoperator}.
       \end{itemize}

 \par In the Banach lattice context, Aqzzouz and Bouras \cite{order DP operator},  El Kaddouri et al. \cite{orderlimitedoperators} and H'michane et al. \cite{w*DPoperator} considered the operators which take order bounded subsets of a Banach lattice to Dunford-Pettis (resp. limited) sets, which are called \textit{order Dunford-Pettis} (resp. \textit{order limited}) \textit{operators}. By the well-known Riesz-Kantorovic formula, we can easily see that a bounded linear operator $T:E\to X$ is order Dundord-Pettis (resp. order limited) if and only if $|T^{\prime}x_n'|\xrightarrow{w^*}0$ in $E'$ whenever $x_n'\xrightarrow{w}0$ (resp. $x_n'\xrightarrow{w^*}0$) in $X'$.

  Note that every relatively weakly compact set and every order bounded set in a Banach lattice belongs to  a more general class of sets, the so-called disjointly weakly compact sets. Following W. Wnuk \cite{WOrder}, we  call a  bounded subset $A$ of a Banach lattice $E$ a \textit{disjointly weakly compact set}   if every disjoint sequence from $Sol(A)$ converges weakly to zero. Many other classes of sets  in a Banach lattice, e.g., weakly precompact sets, Dunford-Pettis sets and limited  sets are  disjointly weakly compact sets (see \cite[Remark 2.4]{disjointlyweaklycompact}). It is well known that, for a Banach lattice $E$, $E'$ has order continuous norm if and only if $B_E$ is disjointly weakly compact.  Xiang, Chen and Li \cite{disjointlyweaklycompact}  extensively investigated disjointly weak compactness  properties in Banach lattices.  Recently, Chen and Li \cite{RdpsetandV*sets} proved that  a bounded subset $K$ of  a Banach lattice $E$ is  disjointly weakly compact if and only if $T(K)$ is a relatively weakly compact set for every Banach space $X$ and every Dunford-Pettis operator $T:E\rightarrow X$, if and only if $K$ is a V$^*$-set of A. Pelczy\'nski. Here, a bounded subset $K$ of a Banach space $X$ is called a \textit{V$^{\,*}$-set} if  $\displaystyle\mathop{\textmd{sup}}_{x\in K}|x_n'(x)|\rightarrow0$ for each weakly unconditionally Cauchy series $\sum_{n}x_n'$ in $X'$ (see \cite{V*set}). Therefore, it follows that a bounded linear operator between two Banach lattices always preserves the disjointly weak compactness of a set \cite[Corollary 0.2]{RdpsetandV*sets}.
	
 In \cite{ChenDW}, the authors  of the present paper defined a  class of operators called \textit{$DW$-compact operators}, which map disjointly weakly compact subsets of a Banach lattice to relatively compact sets. Motivated by the notion of $DW$-compact operators we introduce two new classes of operators: $DW$-DP operators and $DW$-limited operators, which carry disjointly weakly compact sets  to Dunford-Pettis sets and limited sets, respectively. Clearly, every $DW$-compact operator is  $DW$-DP and $DW$-limited, and every $DW$-DP operator (resp. $DW$-limited) operator is weak Dunford-Pettis and order Dunford-Pettis (resp. weak$^*$ Dunford-Pettis and order limited). We show that $DW$-DP (resp. $DW$-limited) operators are precisely those operators which are  wDP  and order Dunford-Pettis (resp.  w$^*$DP  and order limited) operators (Theorem \ref{dw-DP equvilent to wDP and order DP} \& \ref{dw-limited equvilent to w*DP and order limited}).
 Furthermore, the approximation properties of positive $DW$-DP and positive $DW$-limited operators are also given (Theorem \ref{approximation of dw-DP operator} \& \ref{approximation of dw-limited operator}).

	\section{$DW$-DP operators and $DW$-limited operators}\label{main}
	\setcounter{equation}{0}
	\begin{definition}
		A bounded linear operator $T:E\rightarrow X$ from a Banach lattice $E$ to a Banach space $X$ is called a $DW$-DP (resp. $DW$-limited) operator if $T(A)$ is  a Dunford-Pettis (resp. limited) set for each disjointly weakly compact subset $A$ of $E$.
	\end{definition}
	
	Recall that a bounded linear operator $T:X\rightarrow E$ is a \textit{disjointly weakly compact operator} if $T(B_X)$ is a disjointly weakly compact subset of $E$ \cite[Definition 3.1]{disjointlyweaklycompact}. $DW$-DP operators can be characterized by  Dunford-Pettis operators and disjointly weakly compact operators. By saying that $(x_n)_{n=1}^{\infty}\subset E$ is a  disjointly weakly compact sequence, we mean that the set $\{x_n:n\in \mathbb{N}\}$ is disjointly weakly compact.
	\begin{theorem}\label{characterization of dw-DP operators}
		For a bounded linear operator $T:E\rightarrow X$ from a Banach lattice to a Banach space the following statements are equivalent:
		\begin{enumerate}
			\item  $T:E\rightarrow X$ is a $DW$-DP operator.
			\item For each disjointly weakly compact sequence $(x_n)$ of $E$ and each weakly null sequence $(f_n)$ of $X'$, $f_n(Tx_n)\rightarrow0$.
			\item For an arbitrary Banach space $Y$ and each disjointly weakly compact operator $S:Y\rightarrow E$, the adjoint of $TS$ is a Dunford-Pettis operator.
			\item For each disjointly weakly compact operator $S:l^1\rightarrow E$, the adjoint of $TS$ is a Dunford-Pettis operator.
		\end{enumerate}
	\end{theorem}
	\begin{proof}
		(1)$\Rightarrow$(2) Assume that $(x_n)$ is an arbitrary disjointly weakly compact sequence of $E$. Then  $\{T(x_n): n\in N\}$ is a Dunford-Pettis set. Therefore, for each weakly null sequence $(f_n)$ in $X'$, it follows that  $|f_n(Tx_n)|\leq\displaystyle\mathop{\textmd{sup}}_{k}|f_n(Tx_k)|\rightarrow0 \,\,~~~(n\rightarrow0)$.
		
		(2)$\Rightarrow$(3) It suffices to prove that $\|(TS)'(f_n)\|\rightarrow0$ whenever $f_n\xrightarrow{w}0$ in $X'$. Assume by way of contradiction that $\|(TS)'(f_n)\|>\varepsilon_0$ for some weakly null sequence $(f_n)\subset X'$ and  some $\varepsilon_0>0$. Then for each $n\in \mathbb{N}$, there exists some $y_n$ of $B_Y$ such that$$|(TS)'(f_{n})(y_n)|=|f_{n}(TS(y_n))|\geq\varepsilon_0.$$
		Note that  $(S(y_n))$ is a disjointly weakly compact sequence since  $S$ is  disjointly weakly compact. Therefore, $|\langle f_{n},T(Sy_n)\rangle|\rightarrow0~~~(n\rightarrow\infty)$, which leads to a contradiction.
		
		(3)$\Rightarrow$(4) is obvious.
		
		(4)$\Rightarrow$(1) Let $A$ be a disjointly weakly compact subset of $E$ and let $(f_n)$ be an arbitrary weakly null sequence of $X'$. For an arbitrary sequence $(x_n)$ from $A$, let us  define $S:l^1\rightarrow E$ by $S((\alpha_1, \alpha_2, \alpha_3, \ldots))=\sum_{1}^{\infty}\alpha_nx_n$. Then $S$ is a disjointly weakly compact operator. See the proof of \cite[Theorem 2.4]{ChenDW}. Therefore,  $(TS)'$ is a Dunford-Pettis operator, and hence we have $\|(TS)'(f_n)\|\rightarrow0$. Then $$|f_n(Tx_n)|=|(TS)'(f_n)(e_n)|\rightarrow0$$  This implies that  $T(A)$ is a Dunford-Pettis set and hence $T$ is $DW$-DP.
	\end{proof}
	
	In \cite{ChenDW} the authors introduced the (d)-DP  and (d)-DP$^*$ properties. A Banach lattice $E$ is said to have the \textit{(d)-DP} (resp. \textit{(d)-DP$^{\,*}$}) \textit{property} if every disjointly weakly compact subset of $E$ is a Dunford-Pettis (resp. limited) set, or equivalently, if each weakly compact (resp. continuous) operator from $E$ to $c_0$ is $DW$-compact. As an immediate consequence of Theorem \ref{characterization of dw-DP operators} the following is clear.
	
	\begin{corollary}
		For a Banach lattice $E$ the following assertions are equivalent.
		\begin{enumerate}
			\item $E$ has the (d)-DP property, i.e., the identity operator $I:E\rightarrow E$ is a $DW$-DP operator.
			\item For each disjointly weakly compact sequence $(x_n)$ of $E$, each weakly null sequence $(f_n)$ of $E'$, $f_n(x_n)\rightarrow 0$.

			\item For an arbitrary Banach space $Y$ and each disjointly weakly compact operator $S:Y\rightarrow E$, $S'$ is a Dunford-Pettis operator.
			\item For each disjointly weakly compact operator $S:l^1\rightarrow E$,  $S'$ is a Dunford-Pettis operator.
		\end{enumerate}
	\end{corollary}	
	
	It should be noted that a bounded set $A$ of a Banach space $X$ is a Dunford-Pettis set if and only if every weakly compact operator from $X$ into an arbitrary Banach space carries $A$ onto a relatively compact set (see, e.g., \cite[p.350]{positiveoperators}). Next, we can see that  $DW$-DP operators can also be characterized in terms of $DW$-compact operators.
	
	\begin{proposition}\label{characterization of dw-DP used by DW}
		For a bounded linear operator $T:E\rightarrow X$ from a Banach lattice to a Banach space the following statements are equivalent:
		\begin{enumerate}
			\item $T$ is a $DW$-DP operator.
			\item For each weakly compact operator $S$ from $X$ to an arbitrary Banach space $Y$, $ST$ is a $DW$-compact operator.
			\item For each weakly compact operator $S$ from $X$ to $c_0$, $ST$ is a $DW$-compact operator.
		\end{enumerate}
	\end{proposition}
	\begin{proof}
		(1)$\Rightarrow$(2)$\Rightarrow$(3) are obvious.
		
		(3)$\Rightarrow$(1) Let $(x_n)$ be an arbitrary disjointly weakly compact sequence of $E$ and let $(f_n)$ be a weakly null sequence of $X'$. Define $S:X\rightarrow c_0$ by $S(x)=(f_n(x))_n$. Since $S$ is a weakly compact operator (see, e.g., \cite[Theorem 5.26]{positiveoperators}), by our hypothesis, $ST:E\rightarrow c_0$ is $DW$-compact. Therefore $\{ST(x_n)=(f_k(Tx_n))_{k}\}$ is relatively compact in $c_0$, that is to say, $\displaystyle\mathop{\textmd{sup}}_{n}|f_k(Tx_n)|\rightarrow0$ as $k\rightarrow\infty$. Then we have $f_n(Tx_n)\rightarrow0$, i.e., $T$ is a $DW$-DP operator.
	\end{proof}
	
	Recall that a bounded linear operator $T:X\rightarrow Y$ between Banach spaces $X, Y$ is called a \textit{limited operator} if $TB_X$ is a limited subset of $Y$ \cite{Limitedsets}. Clearly, $T:X\rightarrow Y$ is limited if and only if $\|T'y_n'\|\to 0$ whenever $y_n'\xrightarrow{w^*}0$ in $Y'$. A Banach space in which every limited set is relatively compact is called a \textit{Gelfand-Phillips space}. All separable Banach spaces and all WCG-spaces are Gelfand-Phillips spaces.  A $\sigma$-Dedekind complete Banach lattice $E$ is a Gelfand-Phillips space if and only if
$E$ has order continuous norm (see, e.g., \cite[Theorem 4.5]{WOrder}). For the case of $DW$-limited operators, we can obtain the analogues of the preceding results and omit the proofs.
 	
	\begin{theorem}\label{char of dw-limited by sequence}
		For a bounded linear operator $T:E\rightarrow X$ from a Banach lattice to a Banach space the following statements are equivalent:
		\begin{enumerate}
			\item  $T$ is a $DW$-limited operator.
			\item For each disjointly weakly compact sequence $(x_n)$ of $E$ and each weak$^*$-null sequence $(f_n)$ of $X'$, we have $f_n(T(x_n))\rightarrow 0$.
			\item For an arbitrary Banach space $Y$ and each disjointly weakly compact operator $S:Y\rightarrow E$, $TS$ is a limited operator.
			\item For each disjointly weakly compact operator $S:l^1\rightarrow E$, $TS$ is a limited oeprator.
		\end{enumerate}
		\end{theorem}
	
	\begin{corollary}
		For a Banach lattice $E$ the following assertions are equivalent.
		\begin{enumerate}
			\item $E$ has the (d)-DP$^{\,*}$ property, i.e., the identity operator $I:E\rightarrow E$ is a $DW$-limited operator.
			\item For each disjointly weakly compact sequence $(x_n)$ of $E$, each weak$^{\,*}$-null sequence $(f_n)$ of $E'$, $f_n(x_n)\rightarrow 0$.
			\item For an arbitrary Banach space $Y$ and each disjointly weakly compact operator $S:Y\rightarrow E$, $S$ is a limited operator.
			\item For each disjointly weakly compact operator $S:l^1\rightarrow E$, $S$ is a limited operator.
		\end{enumerate}
	\end{corollary}
	
	\begin{proposition}
		For a bounded linear operator $T:E\rightarrow X$ from a Banach lattice to a Banach space the following statements are equivalent:
		\begin{enumerate}
			\item $T$ is a $DW$-limited operator.
			\item For each bounded operator $S$ from $X$ to an arbitrary Gelfand-Phillips  space $Y$, $ST$ is a $DW$-compact operator.
			\item For each bounded operator $S$ from $X$ to $c_0$, $ST$ is a $DW$-compact operator.
		\end{enumerate}
	\end{proposition}

	\begin{remark}\label{remark a}
		(1) Clearly, if $E$ has the (d)-DP (resp. (d)-DP$^*$) property, then every bounded linear operator $T:E\rightarrow X$ from $E$ to an arbitrary Banach space $X$ is a $DW$-DP (resp. $DW$-limited) operator. For two Banach lattices $E$ and $F$, if either $E$ or $F$ has the (d)-DP property (resp. (d)-DP$^*$ property), then every bounded linear operator from  $E$ to $F$ is $DW$-DP (resp. $DW$-limited) since the disjointly weak compactness can be  preserved by bounded linear operators between Banach lattices \cite[Corollary 0.2]{RdpsetandV*sets}.
		
		(2)\label{dw-DP is not dw-limied} It is obvious that every $DW$-limited operator is a $DW$-DP operator, but the converse is  not necessarily true. The identity operator $I:c\rightarrow c$ is an example of a $DW$-DP operator which is not $DW$-limited.
	
    	(3) Clearly, every $DW$-DP (resp. $DW$-limited ) operator from a Banach lattice to a Banach space is a wDP (w$^*$DP) operator. The identity operator $I:l^\infty \rightarrow l^\infty$ is a w$^*$DP operator,  but $I$ is not a $DW$-DP operator since $l^\infty$ has the DP$^*$ property, but $l^\infty$ dose not have the (d)-DP property \cite[Remark 3.5]{ChenDW}. It should be noted that a Banach lattice $E$ is a KB-space if and only if each disjointly weakly compact set of $E$ is relatively weakly compact \cite[Proposition 2.6]{disjointlyweaklycompact}. Therefore,  every wDP (resp. w$^*$DP)  operator from a KB-space to an arbitrary Banach space is $DW$-DP (resp. $DW$-limited).
	\end{remark}
	
Note that  every $DW$-limited operator into $c_0$ is $DW$-compact since $c_0$ is a Gelfand-Phillips space. Therefore, the following characterization of the (d)-DP$^{*}$ property tells us when every $DW$-DP operator is  $DW$-limited.

\begin{proposition}\label{every operator is DW-limited}
For a Banach lattice $E$ the following assertions are equivalent.
		\begin{enumerate}
			\item $E$ has the (d)-DP$^{\,*}$ property.
			\item Each bounded linear operator from $E$ into an arbitrary Banach space $X$ is a $DW$-limited operator.
			\item Each $DW$-DP operator from $E$ into an arbitrary Banach space $X$ is a $DW$-limited operator.
			\item Each $DW$-DP operator from $E$ into $c_0$ is a $DW$-limited (and hence $DW$-compact) operator.
		\end{enumerate}
	\end{proposition}

\begin{proof}
It suffices to prove that $(4)\Rightarrow(1)$. We claim that every bounded linear operator $T:E\to{c_0}$ is $DW$-DP. To this end, let $A$ be a  disjointly weakly compact subset  of $E$. Then $T(A)$ is likewise disjointly weakly compact (see \cite[Corollary 0.2]{RdpsetandV*sets}). Therefore, $T(A)$ is a Dunford-Pettis set since $c_0$ has the (d)-DP property \cite[Remark 3.5 (1)]{ChenDW}. This implies that $T$ is a $DW$-DP operator. Hence, our desired result follows from Proposition 3.2 of \cite{ChenDW}.
	\end{proof}

It is easy to see that  Dunford-Pettis sets  coincide with limited sets in a Grothendieck space. The class of weak Grothendieck operators was introduced by Oughajji  et al. in \cite{wGrothendieck}. An operator $T:X\to Y$  is called a \textit{weak Grothendieck operator} if $f_n(Tx_n)\rightarrow0$ for every weakly null and Dunford-Pettis sequence $(x_n)$ of $X$ and every weak$^*$-null sequence $(f_n)$ of $Y'$, or equivalently, if $T$ carries Dunford-Pettis sets of $X$ onto limited sets of $Y$. Obviously, every $DW$-limited operator is a weak Grothendieck operator since every Dunford-Pettis subset of a Banach lattice is disjointly weakly compact \cite[Remark 2.4 (1)]{disjointlyweaklycompact}. However, the  converse does not hold. The identity operator $I:l^\infty\rightarrow l^\infty$ is an example of a weak Grothendieck operator which is not a $DW$-limited operator since $l^\infty$ is a Grothendieck space without  the (d)-DP$^*$ property.

	\begin{proposition}
For a Banach lattice $E$  the following assertions are equivalent.
		\begin{enumerate}
			\item[{\rm (1)}] Each weak Grothendieck operator from $E$ into an arbitrary Banach space $X$ is a $DW$-limited operator.
			\item[{\rm (2)}] Each weak Grothendieck operator from $E$ into $c_0 $ is a $DW$-limited operator.
			\item[{\rm (3)}] $E$ has the (d)-DP property.

		\end{enumerate}
	\end{proposition}
	\begin{proof}
		It suffices to prove $(2)\Rightarrow(3)$. Assume by way of contradiction that $E$ does not have the (d)-DP property. Then there exists a disjointly weakly compact subset $A$ of $E$ which is not a Dunford-Pettis set. This implies that  there exits a weakly null sequence $(f_n)$ of $E'$ such that $\displaystyle\mathop{\textmd{sup}_{x\in A}}|f_n(x)|\geq\epsilon_0$ for all $n\in \mathbb{N}$ and some $\epsilon_{0}>0$. Let us define $T:E\rightarrow c_0$ by $T(x)=(f_n(x))$ for all $x\in E$. Then $T$ is a weakly compact operator and hence $T$ is a weak Grothendieck operator (see, e.g., \cite[Theorem 5.26, Theorem 5.98]{positiveoperators}). However, $T(A)$ is not relatively compact and hence $T(A)$ is not limited in  $c_0$ since $c_0$ is a Gelfand-Phillips space. This implies that $T$ is not $DW$-limited, which leads to a contradiction.
	\end{proof}
	
	It should be noted that a wDP operator is not necessarily a w$^*$DP operator. For instance, the identity operator $I:c_0\rightarrow c_0$ is  $DW$-DP  and hence it is wDP . However,  it is not a w$^*$DP operator since $c_0$ does not have the DP$^*$ property.
	
	\begin{proposition}\label{wDP and w*DP}
		For a Banach lattice $X$  the following assertions are equivalent.
		\begin{enumerate}
			\item[{\rm (1)}] $X$ has the DP$^{\,*}$ property.
            \item[{\rm (2)}] Each bounded linear operator from $X$ into an arbitrary Banach space  is a w$^*$DP operator.
			\item[{\rm (3)}] Each wDP operator from $X$ into an arbitrary Banach space  is a w$^*$DP operator.
			\item[{\rm (4)}] Each wDP operator from $X$ into $c_0$ is a w$^*$DP operator.
		\end{enumerate}
	\end{proposition}
	\begin{proof}
		It suffices to prove that $(4)\Rightarrow(1)$. Assume by way of contradiction that  there exists a relatively weakly compact subset $A$ of $X$ which is not limited. Then we can find a weak$^*$-null sequence $(f_n)$ of $X'$ such that  $\displaystyle\mathop{\textmd{sup}}_{x\in A}|f_n(x)|>\epsilon_0$ for all $n$ and some $\varepsilon_0>0$. Define $T:X\rightarrow c_0$ by $T(x)=(f_n(x))$. For each relatively weakly compat subset $B$ of $X$, $T(B)$ is likewise relatively weakly compact, and  hence $T(B)$ is a Dunford-Pettis set since $c_0$ has the DP property. This implies that $T$ is a wDP operator. However,  $T(A)$ is not limited in $c_0$. Therefore, $T$ is not a w$^*$DP operator.
	\end{proof}

By definition, every $DW$-DP  (resp. $DW$-limited) operator is wDP and order Dunford-Pettis (resp. w$^*$DP and order limited).	As we have pointed out in Remark \ref{remark a} (3), a w$^*$DP operator  is not necessarily a $DW$-DP operator. Furthermore, an order limited operator need not be $DW$-DP operator. For instance, since $l^p$ $(1<p<\infty)$ is a discrete Banach lattice with order continuous norm and  every order interval in $l^p$ is compact, the identity operator $I:l^{p}\to l^{p}$ is order limited. However, $I:l^{p}\to l^{p}$ is not $DW$-DP (since $B_{l^{p}}$ is disjointly weakly compact, but $B_{l^{p}}$ is not a Dunford-Pettis set.) The following result shows that $DW$-DP operators are precisely those operators which are both wDP and order Dunford-Pettis.
	
	\begin{theorem}\label{dw-DP equvilent to wDP and order DP}
		Let $T:E\rightarrow X$ be a bounded linear operator from a Banach lattice $E$ into a Banach space $X$. Then $T$ is a $DW$-DP operator if and only if $T$ is a wDP  and order Dunford-Pettis operator.
	\end{theorem}
	\begin{proof}
		Assume that  $T:E\rightarrow X$ is a wDP and order Dunford-Pettis operator. Let $A$ be a disjointly weakly compact subset of $E$. By definition, we can assume without loss of generality that $A$ is solid. For every weakly null sequence $(f_n)$ of $X'$, we have to show that $(f_n)$ converges uniformly to zero on $T(A)$. To this end, let us define $\rho(x)=\displaystyle\mathop{\textmd{sup}}_{k}|f_k(x)|$ for every $x\in X$. Then $\rho$ is a norm continuous seminorm on X.
		
	For each disjoint sequence $(x_n)$ of $A$, we have $x_{n}\xrightarrow{w}0$ in $E$ since $A$ is a disjointly weakly compact. We claim that $\rho(Tx_n)\rightarrow0$~~$(n\to\infty)$. Otherwise, by passing to a subsequence if necessary, there exists some $\epsilon_0>0$  such that $\rho(Tx_n)>\epsilon_0$ for each $n$. Then, for each $n$, we can find some $k_n\in \mathbb{N}$ such that $|f_{k_n}(Tx_n)|>\epsilon_0$. We can easily see that the set $\{k_n: n\in N\}$ is not finite since $Tx_{n}\xrightarrow{w}0$. Hence, without loss of generality we can assume that $$k_{1}<k_{2}<\cdots<k_{n}<\cdots$$Then $f_{k_n}\xrightarrow{w}0$,  and $\epsilon_{0}<|f_{k_n}(Tx_n)|\to 0~~(n\to\infty)$ since $T$ is a wDP operator. This is impossible. So we have $\displaystyle\mathop{\lim}_{n}\rho(Tx_n)=0$ for each disjoint sequence $(x_n)$ of $A$.  Therefore,  from \cite[Theorem 4.36]{positiveoperators} it follows that, for each $\epsilon>0$, there exists some $u\in E^+$ lying in the ideal generated by $A$ such that $$\rho(T[(|x|-u)^+])<\frac{\epsilon}{2}$$ holds for all $x\in A$. Therefore, from the identities
\begin{eqnarray*}
x=x^{+}-x^{-}&=&[x^{+}\wedge u+(x^{+}-u)^{+}]- [x^{-}\wedge u+(x^{-}-u)^{+}]\\&=& [x^{+}\wedge u-x^{-}\wedge u]+[(x^{+}-u)^{+}-(x^{-}-u)^{+}]
\end{eqnarray*}
and $x^{+}\wedge u-x^{-}\wedge u\in [-u,u]$, we can see that  $T(A)\subset T[-u,u]+\epsilon B_{\rho}$, where $B_{\rho}=\{z\in X:\rho(z)\leq1\}$.
Note that $T[-u,u]$	is a Dunford-Pettis set since $T$ is an order Dunford-Pettis operator. Hence, $(f_n)$ converges uniformly to zero on $T[-u,u]$.	
 For each $x\in A$, there exists some $y\in [-u,u]$, some $z\in B_{\rho}$ such that $Tx=Ty+\epsilon z$. Therefore, from $$|f_n(Tx)|\leq |f_n(Ty)|+\epsilon\leq \displaystyle\mathop{\textmd{sup}}_{y\in [-u,u]}|f_n(Ty)|+\epsilon$$it follows that $\displaystyle\mathop{\limsup}_{n}\{|f_{n}(Tx)|:x\in A\}\leq\epsilon$. Since $\epsilon>0$ is arbitrary,  $(f_n)$ converges uniformly to zero on $T(A)$. This implies that  $T(A)$ is a Dunford-Pettis set, and hence $T$ is a $DW$-DP operator.
 \end{proof}
	
Similarly, we have the following result and the proof is almost the same as that of Theorem \ref{dw-DP equvilent to wDP and order DP} if   ``$f_{n}\xrightarrow{w}0$" is replaced with ``$f_{n}\xrightarrow{w^*}0$".	
	
\begin{theorem}\label{dw-limited equvilent to w*DP and order limited}
		Let $T:E\rightarrow X$ be a bounded linear operator from a Banach lattice $E$ into a Banach space $X$. Then $T$ is a $DW$-limited operator if and only if $T$ is a w$^{\,*}$DP and order limited operator.
	\end{theorem}

Next, we come to the domination problem of $DW$-DP operators and $DW$-limited operators. Let $S, T:E\to F$ be two positive operators between Banach lattices $E$ and $F$ such that $0\leq S\leq T$.  Kalton and Saab \cite{idealproperties} proved that $S$ is also a wDP operator if  $T$ is a wDP operator. Chen et al. \cite{dominationofpositiveweakstar} proved that if $F$  is $\sigma$-Dedekind complete and $T$ is a  w$^*$DP operator, then $S$ is likewise w$^*$DP operator . For the domination properties of positive order Dunford-Pettis operators and order limited operators on Banach lattices, H'michane and El Fahri \cite{dominationoflimitedandorderDP} have proved that every positive operator dominated by a positive order Dunford-Pettis operator is also order Dunford-Pettis, and if the range space is $\sigma$-Dedekind complete, then every positive operator dominated by a positive order limited operator is also order limited. Then, As immediate consequences of  Theorem \ref{dw-DP equvilent to wDP and order DP} and Theorem \ref{dw-limited equvilent to w*DP and order limited}, we have the following results about the domination properties of  $DW$-DP operators and $DW$-limited operators.
	
	\begin{corollary}
		Let $S$ and  $T$ be two positive operators between two Banach lattices $E$ and  $F$ such that  $0\leq S\leq T$. If $T$ is a $DW$-DP operator, then $S$ is likewise a $DW$-DP operator.
	\end{corollary}
	
	\begin{corollary}\label{domination of dw-limited}
		Let $E$ and $F$ be two Banach lattices with $F$ $\sigma$-Dedekind complete and  let $S, T:E\to F$ be two positive operators  such that $0\leq S\leq T$. If $T$ is a $DW$-limited operator, then $S$ is also a $DW$-limited operator.
	\end{corollary}

The following example shows that the condition ``\,$F$ is $\sigma$-Dedekind complete\," is essential in Corollary \ref{domination of dw-limited}.
\begin{example}
Let us consider the case when $E=L^{1}[0,1]$ and $F=c$. Since $L^{1}[0,1]$ is a $KB$-space, disjointly weakly compact sets and relatively weakly compact sets coincide in $L^{1}[0,1]$ \cite[Proposition 2.6]{disjointlyweaklycompact}. Therefore, the class of $DW$-limited operators and the class of  Dunford-Pettis operators from $L^{1}[0,1]$ into $c$ are the same (since $c$ is a Gelfand-Phillips space). Now, since $L^{1}[0,1]$ does not have weakly sequentially continuous lattice operations and  $c$ does not have  order continuous norm,  from the converse for the Kalton-Saab theorem proved by W. Wickstead \cite{Wickstead} it follows that there exists a positive operator $S:L_{1}[0,\,1] \to c$ such that $S$ is dominated by a Dunford-Pettis operator while $S$ is not Dunford-Pettis. This implies that the assertion in Corollary \ref{domination of dw-limited} does not necessarily hold when the target space is not $\sigma$-Dedekind complete.
\end{example}
\vskip 5mm

	\section{The approximation properties of  $DW$-DP operators and  $DW$-limited operators}

		Aliprantis and Burkinshaw gave a result on the approximation property of  positive wDP operators (see, e.g., \cite[Theorem 5.100]{positiveoperators}).  For the approximation property of  positive w$^*$DP operators see  Theorem 2.3 \& 2.4 of Chen et al. \cite{dominationofpositiveweakstar}. In this section, we consider the approximation properties of  positive $DW$-DP operators and  positive $DW$-limited operators.
\begin{theorem}\label{approximation of dw-DP operator}

		Let $T:E\rightarrow F$ be a positive $DW$-DP operator between Banach lattices. If $W\subset E$  and $V\subset F'$ are two  disjointly weakly compact sets, then the following statements hold:
		\begin{enumerate}
			\item For every disjoint sequence $(x_n)$ in the solid hull of $W$, the sequence $(Tx_n)$ converges uniformly to zero on the solid hull of $V$.
			\item For each $\epsilon>0$ there exists some $u\in E^+$ satisfying
			\begin{center}
				$\langle|f|,T[(|x|-u)^+]\rangle<\epsilon$
			\end{center}
			for all $x\in W$ and all $f\in V$.
		\end{enumerate}
	\end{theorem}

\begin{proof}
 By definition, we can assume  without loss of generality that $W$ and $V$ are solid.

 (1) Let $(x_n)$ be a disjoint sequence in $W\cap E^+$. Then $0\leq x_{n}\xrightarrow{w}0$ in $E$ since $W$ is disjointly weakly compact, and hence the positivity of $T$ implies that $0\leq Tx_{n}\xrightarrow{w}0$ in $F$. Assume by way of contradiction that $\displaystyle\mathop{\sup}_{f\in V}|f(Tx_n)|\nrightarrow0$. By passing to a subsequence if necessary, there would exist a  sequence $(f_n)\subset V$ satisfying ${\vert{f_n}\vert}(Tx_n)\geq\vert f_{n}(Tx_n)\vert>{\epsilon}$  for some $\epsilon>0$ and for all $n\in\mathbb{N}$. Since $0\leq Tx_n\stackrel{w}{\rightarrow}0$, by induction there exists a strictly increasing subsequence $(n_{k})_{k=1}^{\infty}$ of $\mathbb{N}$ such
that $\langle4^{m}\sum_{k=1}^{\,m}\vert f_{n_{k}}\vert,  Tx_{n_{m+1}}\rangle<2^{-m}$ for all $m\in\mathbb{N}$. See the proof of \cite[Theorem 2.5]{almostlimitedsets} for details.  Let $$f=\sum_{k=1}^{\,\infty}2^{-k}\vert f_{n_{k}}\vert,\quad g_{m}=\left(\vert f_{n_{m+1}}\vert-4^{m}\sum_{k=1}^{\,m}\,\vert f_{n_{k}}\vert-2^{-m}f\right)^{+}.$$Then,
by Lemma 4.35 of \cite{positiveoperators} $(g_{m})$ is a disjoint sequence in $V$ since $V$ is solid. Now we have
\begin{eqnarray*}
    g_{m}(Tx_{n_{m+1}})&=&\left\langle\left(|f_{n_{m+1}}|-4^{m}\sum_{k=1}^{\,m}\,|f_{n_{k}}|-2^{-m}f\right)^{+}, \,Tx_{n_{m+1}}\right\rangle\\&\geq& \left\langle\left(|f_{n_{m+1}}|-4^{m}\sum_{k=1}^{\,m}\,|f_{n_{k}}|-2^{-m}f\right),\, Tx_{n_{m+1}}\right\rangle
    \\&=&|f_{n_{m+1}}|(Tx_{n_{m+1}})-\left\langle4^{m}\sum_{k=1}^{\,m}\,|f_{n_{k}}|, Tx_{n_{m+1}}\right\rangle-2^{-m}f(Tx_{n_{m+1}})\\&>&\epsilon-2^{-m}-2^{-m}
    f(Tx_{n_{m+1}}).
\end{eqnarray*}
 Therefore,  $g_{m}(Tx_{n_{m+1}})>\frac{\epsilon}{2}$  holds when $m$ is sufficiently large.  On the other hand, since  $(g_{m})\subset V$  is a  disjoint sequence and $V$ is disjointly weakly compact, it follows  that $g_m\xrightarrow{w}0$ in $F'$. Then
 $$\frac{\epsilon}{2}<g_{m}(Tx_{n_{m+1}})\leq \sup_{x\in W}|g_{m}(Tx)|\rightarrow0\,\,~~~~\,\, (n\to\infty)$$since $T(W)$ is a Dunford-Pettis set. This leads to a contradiction. Therefore, $(Tx_n)$ converges uniformly to zero on  $V$.
 \par Now, for the general case, let $(x_n)$ be a disjoint sequence in $W$. Then $(x_{n}^+)$ and $(x_{n}^-)$ are disjoint sequences in $W$. Therefore, $(Tx_{n}^+)$ and $(Tx_{n}^-)$ converge uniformly to zero on $V$. It follows that $Tx_n$ converges uniformly to zero on $V$, as desired.

 \par(2) Let us define a seminorm $\rho_{V}$ on $F$ as $\rho_{V}(y)=\displaystyle\mathop{\sup}_{f\in V}|f(y)|$ for $y\in F$. Since $V$ is solid, by the Riesz-Kantorovich Formula we have $\rho_{V}(y)=\displaystyle\mathop{\sup}_{f\in V} |f|(|y|)$. It is obvious that $\rho_{V}$ is a norm continuous seminorm on $F$. From Part (1) and \cite[Theorem 4.36]{positiveoperators} it follows that for each $\epsilon>0$ there exists some $u\in E^+$ lying in the ideal generated by $W$  such that $\rho_{V}\left( T[(|x|-u)^+] \right)<\epsilon$ for all $x\in W$, or equivalently, $T(W)\subset [-Tu,Tu]+\epsilon B_{\rho_{V}}$. Therefore,
 \begin{center}
				$\langle|f|,T[(|x|-u)^+]\rangle<\epsilon$
			\end{center}
			for all $x\in W$ and all $f\in V$.
\end{proof}

 Let us recall that a bounded linear operator $T:E\to X$ from a Banach lattice into a Banach space is called an
\textit{almost Dunford-Pettis operator} if $\|Tx_n\|\to0$ for each disjoint sequence $(x_n)$ in $E$ satisfying $x_{n}\xrightarrow{w}0$.  Chen and Li \cite{V*set} proved that an almost Dunford-Pettis operator carries disjointly weakly compact sets onto relatively weakly compact sets.  The identity operator $I:c_0\rightarrow c_0$ is  an example of a $DW$-DP operator which  is not almost Dunford-Pettis.
\begin{corollary}
Let $E$ and $F$ be two Banach lattices such that $F''$ has order continuous norm. Then every positive $DW$-DP operator from $E$ into $F$ is an almost Dunford-Pettis operator.
\end{corollary}

\begin{proof}
Let $T:E\rightarrow F$ be a positive $DW$-DP operator and let $(x_n)$ be a disjoint weakly null sequence of $E$. Then, the set $W=\{x_n:n\in \mathbb{N}\}$ is  relatively weakly compact and hence disjointly weakly compact. Since $F''$ has order continuous norm, $V=B_{F'}$ is a disjointly weakly compact set. From Theorem \ref{approximation of dw-DP operator} it follows that $(Tx_n)$ converges uniformly to zero on $V$, that is, $\|Tx_n\|\to 0$.
\end{proof}

For our convenience, we call a bounded subset $A$ of $E^\prime$ \textit{disjointly weak$^{\,*}$-compact} if $x_{n}^{\prime}\xrightarrow{w^*}0$ for every disjoint sequence  $(x_{n}^{\prime})\subset Sol(A)$. We know that a relatively weakly compact subset of a  Banach lattice $E$ is  disjointly weakly compact \cite[Theorem 4.34]{positiveoperators}. However, a bounded subset of $E'$ is not necessarily disjointly weak$^{\,*}$-compact although it is certainly relatively weak$^*$-compact by Alaoglu's theorem. It is well known that a Banach lattice $E$ has order continuous norm if and only if $B_{E^{\prime}}$ is disjointly weak$^{*}$-compact (see, e.g., \cite[Corollary 2.4.3]{Meyer}).	
	
\begin{theorem}\label{approximation of dw-limited operator}

		Let $T:E\rightarrow F$ be a positive $DW$-limited operator between Banach lattices. If $W\subset E$  is a disjointly weakly compact set and $V\subset F'$ is a  disjointly weak$^{\,*}$-compact set, then the following statements hold:
		\begin{enumerate}
			\item For every disjoint sequence $(x_n)$ in the solid hull of $W$, the sequence $(Tx_n)$ converges uniformly to zero on the solid hull of $V$.
			\item For each $\epsilon>0$ there exists some $u\in E^+$ satisfying
			\begin{center}
				$\langle|f|,T[(|x|-u)^+]\rangle<\epsilon$
			\end{center}
			for all $x\in W$ and all $f\in V$.
		\end{enumerate}
	\end{theorem}

\begin{proof}
The proof is almost the same as that of Theorem \ref{approximation of dw-DP operator}. We follow the notation used in the proof of Theorem \ref{approximation of dw-DP operator}. Here, the disjointly weak$^*$-compactness of $V$ guarantees that the disjoint sequence  $(g_m)\subset V$ constructed in the proof of Theorem \ref{approximation of dw-DP operator} is weak$^*$-null. Since $T(W)$ is a limited set in $F$, we have $\sup_{x\in W}|g_{m}(Tx)|\rightarrow0\,\,~~~~\,\, (n\to\infty).$
\end{proof}

\end{document}